\documentclass[11pt]{amsart}

\usepackage{amsfonts}
\usepackage{amsmath}
\usepackage{amssymb}
\usepackage{latexsym}
\usepackage{amscd}
\usepackage[mathscr]{eucal}
\usepackage{euler}

\usepackage{youngtab}
\usepackage{gastex}

\numberwithin{equation}{section}

\newtheorem{theorem}{Theorem}[section]

\newtheorem{definition}[theorem]{Definition}
\newtheorem{example}[theorem]{Example}
\newtheorem{lemma}[theorem]{Lemma}

\newtheorem{proposition}[theorem]{Proposition}

\begin{document}
\title{Row convex tableaux and Bott-Samelson varieties}
\author{Philip Foth}
\address{Champlain St. Lawrence College, Quebec, Canada G1V 4K2}
\email{phfoth@gmail.com}

\author{Sangjib Kim}
\thanks{This work was supported by the Ewha Womans University Research Grant of 2013 (Kim)}
\address{Department of Mathematics, Ewha Womans University, Seoul 151-892, Korea}
\email{sk23@ewha.ac.kr}

\begin{abstract}
By using row convex tableaux, we study the section rings of 
Bott-Samelson varieties of type A. We obtain flat deformations 
and standard monomial type bases of the section rings. 
In a separate section, we investigate a three dimensional 
Bott-Samelson variety in detail and compute its Hilbert polynomial 
and toric degenerations.
\end{abstract}

\maketitle


\section{Introduction}


\subsection{}\label{Bquot}

Let $G=\rm{GL}_{n}(\mathbb{C})$ be the general linear
group over the complex number field $\mathbb{C}$ and $B$ be its
Borel subgroup consisting of upper triangular matrices.
For a word $\mathbf{i}=(i_{1},... ,i_{{\ell}})$ with 
$1\leq i_{j}\leq n-1$, the \textit{Bott-Samelson variety} $Z_{\mathbf{i}}$ 
can be defined as the quotient space%
\begin{equation*}
P_{i_{1}}\times P_{i_{2}}\times \cdots \times P_{i_{{\ell}}}/B^{{\ell}}.
\end{equation*}%
Here, $P_{i_{j}}$ is the minimal parabolic subgroup of $G$
associated to the simple reflection 
\begin{equation*}
s_{i_{j}}=(i_{j},i_{j}+1)
\end{equation*}
and $(b_{1},... ,b_{{\ell}}) \in B^{\ell}$ acts on the product of $P_{i_{j}}$'s by
\begin{equation*}
(p_{1},... ,p_{{\ell}}) \cdot (b_{1}, ... , b_{{\ell}})
=(p_{1}b_{1},b_{1}^{-1}p_{2}b_{2}, ... , b_{{\ell}-1}^{-1}p_{{\ell}}b_{{\ell}}).
\end{equation*}%

The Bott-Samelson varieties are defined in \cite{BS55, BS58} and
\cite{De74} to desingularize the Schubert varieties in the flag
manifold $G/B$, and then used to study the Chow ring of $G/B$. In
representation theory, the Bott-Samelson varieties provide Demazure's 
character formula, which can be understood as a generalized Weyl's character formula,
through the section spaces of their line bundles. 

\subsection{}

One can also realize the Bott-Samelson variety $Z_{\mathbf{i}}$ as 
a configuration variety in the product of the Grassmann varieties 
${\rm Gr}(i,n)$ via the map
\begin{align*}
Z_{\mathbf{i}} & \longrightarrow {\rm Gr}(i_{1},n)\times \cdots \times {\rm Gr}(i_{{\ell}},n)\\
(p_1, ..., p_{\ell}) & \longmapsto (p_1 E^{i_1}, p_1 p_2 E^{i_2}, ..., p_1\cdots p_{\ell}E^{i_{\ell}})
\end{align*}
where $E^{i}$ is the $i$-dimensional subspace of $\mathbb{C}^n$ spanned by 
the first $i$ elementary basis elements $\{e_1,..., e_i\}$.
From such a realization, Lakshmibai and Magyar investigated generalized 
Demazure modules and described their standard monomial bases in terms 
of root operators \cite{LM98,Ma98}. See also \cite{LLM02}.

We note that there is a natural line bundle induced from 
the Pl\"{u}cker bundles on the factors ${\rm Gr}(i_{j},n)$, and as is the case 
for the Grassmann varieties and the flag
varieties, we can investigate the Pl\"{u}cker coordinates in terms
of minors over a matrix or Young tableaux and straightening
relations among them. 

In this paper, using the language of row convex tableaux
introduced by Taylor \cite{Ta01}, we study the section rings of 
the Bott-Samelson varieties and its explicit standard monomial type 
bases which are different from the ones given in \cite{LM98,Ma98}.
For $\mathbf{i}$ in \eqref{reduced word} and 
$\mathbf{m}=(m_{1},\cdots ,m_{{\ell}})\in \mathbb{Z}_{\geq 0}^{{\ell}}$, 
our main results are as follows.

\begin{theorem}
Let $\mathsf{M}(\mathbf{i}, d\mathbf{m})$ be the space spanned by
tableaux of shape $(\mathbf{i}, d\mathbf{m})$. The section ring
of the Bott-Samelson variety with respect to the line bundle $L_{\mathbf{m}}$ is
\begin{equation*}
\mathcal{R}_{\mathbf{i,m}} \cong  \bigoplus_{d\geq 0}\mathsf{M}(\mathbf{i}, d\mathbf{m})
\end{equation*}
and straight tableaux of shape $(\mathbf{i}, d\mathbf{m})$ form 
a $\mathbb{C}$-basis of the space $\mathsf{M}(\mathbf{i}, d\mathbf{m})$.
\end{theorem}

Then from SAGBI-Gr\"{o}bner degeneration
techniques (e.g., \cite{MS05, St95}), we obtain a flat degeneration 
of the section ring.

\begin{theorem}
The section ring $\mathcal{R}_{\mathbf{i,m}}$ of the Bott-Samelson variety 
$Z_{\mathbf{i}}$ is a flat deformation of an affine semigroup ring.
\end{theorem}

In the last section, we provide a detailed study of an example for 
the case of $\rm{GL}_{3}(\mathbb{C})$, including toric degenerations, the
corresponding moment polytopes, and computations of the Hilbert
polynomial.
\begin{proposition}
The Hilbert polynomial of the Bott-Samelson variety $Z$ is
\begin{equation*}
{\rm HP}_Z(s) = \frac{5s^3+11s^2+8s+2}{2}.
\end{equation*}
\end{proposition}

In \cite{GK94}, Grossberg and Karshon studied a family of complex
structures on a Bott-Samelson manifold, such that the underlying
real manifold remains the same, but the limit complex manifold
admits a complete, full-dimensional torus action. (They call such
varieties ``Bott towers''. An algebraic version of their construction appeared in
\cite{Pa10}.) Our deformation is algebraic in nature, yet is different, 
as can be seen in examples and also from the fact that in the limit, 
the relationship between ${Z_{\mathbf i}}$ and $G/B$ naturally extends 
to the whole flat family.

\subsection{}

This paper is arranged as follows: in Section \ref{S-section},
we fix notation and some basic definitions which will be useful in the
sequel. Then, we describe the section ring of the Bott-Samelson variety 
in terms of row convex tableaux. In Section \ref{S-flat}, by using the 
fact that straight tableaux form bases of the space of sections, we show 
that the section ring is a flat deformation of an affine semigroup ring. 
In Section \ref{S-straight}, we further investigate straight tableaux
and study their properties. In Section \ref{S-examp}, for 
a three-dimensional Bott-Samelson variety, we compute its toric 
degenerations and Hilbert polynomial.


\section{Row Convex Tableaux and the Section Rings}\label{S-section}


In this section, after introducing row convex tableaux and related notation, 
we describe the section ring of the Bott-Samelson variety associated with
a reduced expression of the longest element of the symmetric group.

\subsection{Row convex tableaux}

A \textit{shape} is a finite collection of pairs of positive integers.
A \textit{tableau} $t$ of shape $D$ is an assignment
of positive integers to elements in $D$:
\begin{equation*}
t: D \longrightarrow \mathbb{Z}_{>0}.
\end{equation*}

One can identify a shape $D$ with a collection of cells arranged in rows and columns
in such a way that there is a cell in the $i$th row and $j$th column 
if and only if $(i,j) \in D$. In this realization, a tableaux of a shape $D$
is a filling of cells in $D$ with positive integers.

\begin{definition}
A \textit{row convex shape} is a shape without gaps in any row. 
That is, if $(r,i)$ and $(r,k)$ are in a shape $D$, then 
$(r,j) \in D$ for all $i<j<k$.
A \textit{row convex tableaux} is a filling of 
a row convex shape with positive integers.
\end{definition}

All the row convex shapes in this paper satisfy the following conditions: the higher rows end at least 
as far to the right as lower rows. Such shapes may be understood as a generalization of 
skew Young diagrams in the following sense (cf. \cite{Ta01}). For two Young diagrams 
\begin{align*}
\lambda=(\lambda_1, ..., \lambda_\ell) \in \mathbb{Z}^{\ell} \text{\ \ such that\ }  
\lambda_1 \geq \cdots \geq \lambda_\ell \geq 0, \\
\mu=(\mu_1, ..., \mu_{\ell}) \in \mathbb{Z}^{\ell} \text{\ \ such that\ }  
\mu_1 \geq \cdots \geq \mu_{\ell} \geq 0
\end{align*}
with $\lambda_i \geq \mu_i$ for all $i$, 
a skew Young diagram $\lambda/\mu$ is the set-theoretic difference 
of the Young diagrams of $\lambda$ and $\mu$. 
If we replace a Young diagram $\mu$ with a sequence of non-negative integers 
$m=(m_1,...,m_{\ell})$ with $\lambda_i \geq m_i$ for all $i$,
then we can obtain a row convex shape $\lambda/m$ by removing the first $m_i$ boxes 
in the $i$th row of the Young diagram $\lambda$ for all $i$.

\subsection{Column sets}

Let us consider the following reduced decomposition of the longest
element $w_0$ in
$\mathfrak{S}_{n}$:%
\begin{equation*}
\underline{w}_0=\underline{w}^{(n)}_0=(s_{1})(s_{2}s_{1})(s_{3}s_{2}s_{1})\cdots (s_{n-1}s_{n-2}\cdots s_{1})
\end{equation*}%
where $s_{i_{j}}$ is the simple reflection $(i_{j},{i_{j}}+1)$.
Note that the length of $\underline{w}_0$ is ${\ell}=n(n-1)/2$. 
Once and for all, we fix the word
\begin{equation}\label{reduced word}
\mathbf{i}=(i_{1},\cdots ,i_{{\ell}})=(1,2,1,3,2,1,\cdots,n-1,n-2,\cdots,1) 
\end{equation}%
associated to the reduced expression $\underline{w}_0=s_{i_{1}}s_{i_{2}}\cdots
s_{i_{{\ell}}}$ of the longest element given above. 

\begin{definition}\label{column sets}
For the reduced word $\mathbf{i}$, the column sets are
\begin{equation*}
C^{(k)}=s_{i_{1}}s_{i_{2}}\cdots s_{i_{k}} [i_k] 
\end{equation*}%
where $[i_k]$ is the set of positive integers not more than $i_k$ 
for $1\leq k\leq {\ell}$.
\end{definition}

Column sets can be defined for any word, but for the reduced word $\mathbf{i}$, 
we can explicitly describe all the column sets. In particular, it is straightforward
to prove that each column set contains consecutive integers.

\begin{lemma} \label{row-convex}
For each $k$, if $a<c<b$ and both $a$ and $b$ are in $C^{(k)}$, then $c\in C^{(k)}$.
To be more precise, for each $j$ with $2\leq j\leq n-1$, let $p_{j}=j(j-1)/2$. 
Then, the column sets are
\begin{equation*}
C^{(p_{j}+t)}=\{t+1,t+2,\cdots ,j+1\}
\end{equation*}%
for $1\leq t\leq j$, and $C^{(1)}=\{2\}$.
\end{lemma}

This shows in particular that if we stack $C^{(k+1)}$ on top of $C^{(k)}$, 
the column sets we defined form a row convex shape
\begin{equation*}
D= \bigcup_{1 \leq k \leq \ell} \left\{(k,c) | c \in C^{(k)} \right\}
\end{equation*}
and its higher rows end at least as far to the right as lower rows.

\begin{example}\label{exii}
For $n=3$, the reduced word is $\mathbf{i}=(121)$ and the column sets are
\begin{align*}
C^{(1)}&=s_1 \{1\}=\{2\}, \\
C^{(2)}&=s_1 s_2 \{1,2\}=\{2,3\}, \\
C^{(3)}&=s_1 s_2 s_1 \{ 1\}= \{3\}. 
\end{align*}
For $n=4$, the reduced word is $\mathbf{i}=(121321)$ and we have
three additional column sets 
\begin{align*}
C^{(4)}&=s_1 s_2 s_1 s_3 \{1,2,3 \}= \{2,3,4\}, \\
C^{(5)}&=s_1 s_2 s_1 s_3 s_2 \{1,2 \} =\{3,4\}, \\
C^{(6)}&=s_1 s_2 s_1 s_3 s_2 s_1 \{ 1\} =\{4\}.
\end{align*}
Then the corresponding row convex shapes for $n=3$ and $n=4$ indicated
by $X$ are respectively
\begin{equation*}
\young(\ \ X,\ XX,\ X\ ) \text{\ \ \  and  \ \ }
\young(\ \ \ X,\ \ XX,\ XXX,\ \ X\ ,\ XX\ ,\ X\ \ ).
\end{equation*}
\end{example}

\subsection{Bott-Samelson varieties}

We will use the realization of the Bott-Samelson variety  
as the variety of configurations of subspaces of $\mathbb{C}^n$. 
For various constructions of the Bott-Samelson varieties and their 
equivalences, we refer the readers to \cite[\S 1]{Ma98}.

For the word $\mathbf{i}$ in \eqref{reduced word}, 
let us write ${\rm Gr}(\mathbf{i})$ for 
\begin{equation*}
{\rm Gr}(i_{1},n)\times \cdots \times {\rm Gr}(i_{{\ell}},n)
\end{equation*}
where ${\rm Gr}(i_k,n)$ is the Grassmann variety of $i_k$ dimensional 
subspaces in $\mathbb{C}^{n}$.

\begin{definition}
The Bott-Samelson variety $Z_{\mathbf{i}}$ is the closure of 
the $B$-orbit of $\underline{x}_{\mathbf{i}}$ in ${\rm Gr}(\mathbf{i})$:
\begin{equation*}
Z_{\mathbf{i}} = \overline{B \cdot \underline{x}_{\mathbf{i}}} 
\subset {\rm Gr}(\mathbf{i})
\end{equation*}%
where $\underline{x}_{\mathbf{i}}=(x_{i_1},..., x_{i_{\ell}})$ is
the point in ${\rm Gr}(\mathbf{i})$ 
whose $k$th coordinate $x_{i_k}$ is the $|C^{(k)}|$-dimensional subspace 
of $\mathbb{C}^{n}$ spanned by the elementary basis elements 
$e_j$ for all $j$ in the column set $C^{(k)}$ for $1 \leq k \leq \ell$.
\end{definition}

We observe that there is a natural line bundle induced from the Pl\"{u}cker bundles 
$\mathcal{O}(1)$ on the factors of ${\rm Gr}(\mathbf{i})$. That is, 
for $\mathbf{m}=(m_{1},\cdots ,m_{{\ell}})\in \mathbb{Z}_{\geq 0}^{{\ell}}$,
we take the powers of the Pl\"{u}cker bundles to obtain an effective line bundle 
on ${\rm Gr}(\mathbf{i})$:
\begin{equation*}
\mathcal{O}(\mathbf{m})= 
\mathcal{O}^{\otimes m_{1}} \otimes \cdots \otimes \mathcal{O}^{\otimes m_{\ell}}.
\end{equation*}

We define the line bundle $L_{\mathbf{m}}$ on the Bott-Samelson variety $Z_{\mathbf{i}}$ 
as the restriction of $\mathcal{O}(\mathbf{m})$ to $Z_{\mathbf{i}} \subset {\rm Gr}(\mathbf{i})$:
\begin{equation}\label{lbundle}
L_{\mathbf{m}}=\mathcal{O}(\mathbf{m})_{| Z_{\mathbf{i}}}.
\end{equation}
and then study the section ring:
\begin{equation*}
\mathcal{R}_{\mathbf{i,m}}=\bigoplus_{d\geq 0}
H^{0}(Z_{\mathbf{i}},L_{\mathbf{m}}^{d}).
\end{equation*}

\subsection{Minors and tableaux}

Let $M_{n}=M_{n}(\mathbb{C})$ be the space of complex $n\times n$ matrices
and $B_{n}=\overline{B}$ be the subspace consisting of upper triangular
matrices:%
\begin{equation*}
B_{n}=\{(x_{ij})\in M_{n}:x_{ij}=0\text{ for }i>j\}.
\end{equation*}%
For $k\leq n$, consider subsets $R=\{r_{1},\cdots ,r_{k}\}$ and 
$C=\{c_{1}, \cdots , c_{k}\}$ of $\{1,\cdots ,n\}$ such that 
$r_{1}<\cdots <r_{k}$ and $c_{1}<\cdots <c_{k}$.
Then, we let $[R:C]$ denote the map from $B_{n}$ to 
$\mathbb{C}$ by assigning to a matrix $b\in B_{n}$ the determinant of the 
$k\times k$ minor of $b$ formed by taking rows $R$ and columns $C$:%
\begin{equation*}
[R:C] =\det \left[
\begin{array}{ccc}
x_{r_{1}c_{1}} & \cdots & x_{r_{1}c_{k}} \\
\vdots & \ddots & \vdots \\
x_{r_{k}c_{1}} & \cdots & x_{r_{k}c_{k}}%
\end{array}%
\right]
\end{equation*}%
where $x_{rc}=0$ if $r>c$.

For subsets $S$ and $S^{\prime}$ of $\{1,\cdots ,n\}$ with the same size, we
can impose a partial ordering: $S\preceq S^{\prime}$ if for each $k$, the $k$%
th smallest element of $S$ is less than or equal to the $k$th smallest
element of $S^{\prime}$. Then note that $[R:C]$ is non-zero only if $%
R\preceq C$. This property is called \textit{flagged}. Since we consider only
minors defined on $B$, from now on we continue to assume this property.

\medskip

By using a Young diagram with a single row consisting of $n$ boxes,
we can record $[R:C]$ by filling in the $c_{i}$th box counting from
left to right with $r_{i}$ for each $i$.
For example, for $n=6$ if $R_{1}=\{1,3,4\}$ and $C_{1}=\{2,3,4\}$ then 
$[R_{1}:C_{1}]$ can be drawn as%
\begin{equation*}
\young(\ 134\ \ )
\end{equation*}%
The product of $k$ of these row tableaux $[R_{i}:C_{i}]$ can be encoded in  
a $k\times n$ rectangular array whose $i$th row counting from bottom to top
is $[R_{i}:C_{i}]$ for $1\leq i\leq k$. 
For example, if $R_{2}=\{2,3,5\},C_{2}=\{3,4,5\},R_{3}=\{4,5\}$ and 
$C_{3}=\{5,6\}$, then $\prod_{1\leq i\leq 3}[R_{i}:C_{i}]$ can be drawn as
\begin{equation*}
\young(\ \ \ \ 45,\ \ 235\ ,\ 134\ \ )
\end{equation*}

\medskip

Next, for ${\ell}=n(n-1)/2$ and
$\mathbf{m}=(m_{1},\cdots ,m_{{\ell}})\in \mathbb{Z}_{\geq 0}^{{\ell}}$, 
consider a collection
\begin{equation*}
\bigcup_{1 \leq k \leq \ell}
\left\{[R_{j}^{(k)}:C^{(k)}]| \ 1\leq j \leq m_{k} \right\}
\end{equation*}%
where $R_{j}^{(k)}$'s are subsets of $\{1,2,...,n \}$ and
$C^{(k)}$'s are the column sets with respect to $\mathbf{i}$ 
(Definition \ref{column sets}). 
Write $|\mathbf{m}|$ for $\sum_{k} m_{k}$. 
Then, by repeating $C^{(k)}$ $m_{k}$ times for
each $k$, the product $\mathsf{t}$ of $[R_{j}^{(k)}:C^{(k)}]$'s can be 
encoded in a $|\mathbf{m}|\times n$ rectangular array having $[R_{j}^{(i)}:C^{(i)}]$ 
as its $(m_{1}+\cdots +m_{i-1}+j)$th row counting from bottom to top.
In this way, we can identify tableaux and products of minors.

\begin{definition}
A \textit{tableaux} $\mathsf{t}$ of shape $(\mathbf{i}, \mathbf{m})$ is
\begin{equation}\label{ttt}
\mathsf{t}=\left( \prod_{1\leq j\leq m_{1}}[R_{j}^{(1)}:C^{(1)}]\right)
\cdot \left( \prod_{1\leq j\leq m_{2}}[R_{j}^{(2)}:C^{(2)}]\right) \cdot
...\cdot \left( \prod_{1\leq j\leq m_{{\ell}}}[R_{j}^{({\ell})}:C^{({\ell})}]\right).
\end{equation} 
\end{definition}

Note that up to sign, we can always assume that 
the entries in each row of $\mathsf{t}$ are increasing from left to right. 
If such is the case, then $\mathsf{t}$ is called a \textit{row standard tableau}.

\medskip

\subsection{Section ring}

From the realization of $Z_{\mathbf{i}}$ as a configuration space in $Gr(\mathbf{i})$, 
we can obtain an explicit description of the space of sections 
$H^{0}(Z_{\mathbf{i}},L_{\mathbf{m}})$ of the line bundle $L_{\mathbf{m}}$. 
In fact, such spaces can be described in a general setting. 
See \cite{LLM02,LM98} and \cite[\S 3]{Ma98} for this direction.

\begin{theorem}\label{coordinate ring}
For $\mathbf{m}=(m_{1},\cdots ,m_{{\ell}})\in \mathbb{Z}_{\geq 0}^{{\ell}}$, 
let $\mathsf{M}(\mathbf{i}, \mathbf{m})$ be the space spanned by tableaux 
of shape $(\mathbf{i}, \mathbf{m})$. Then, we have
\begin{equation*}
\mathsf{M}(\mathbf{i},\mathbf{m}) \cong H^{0}(Z_{\mathbf{i}},L_{\mathbf{m}}).
\end{equation*}
\end{theorem}
\begin{proof}
In the setting
\begin{equation*}
Z_{\mathbf{i}} = \overline{B \cdot \underline{x}_{\mathbf{i}}}
\subset {\rm Gr}(i_{1},n)\times \cdots \times {\rm Gr}(i_{{\ell}},n),
\end{equation*}
the sections of the line bundle $\mathcal{O}(1)$ over the Grassmannian 
${\rm Gr}(i_{k},n)$ can be identified with the maximal minors $\delta_j^{(k)}$ 
defined on the space $X_{k}$ of $n\times i_k$ complex matrices. Therefore, the space 
of sections of $\mathcal{O}(\mathbf{m})$ over ${\rm Gr}(\mathbf{i})$ is 
spanned by the products 
\begin{equation*}
\prod_{j=1}^{m_1} \delta_j^{(1)} \cdot \prod_{j=1}^{m_2} \delta_j^{(2)} 
\cdot ... \cdot \prod_{j=1}^{m_{\ell}} \delta_j^{({\ell})}.
\end{equation*}
We can restrict these sections to $Z_{\mathbf{i}}$ to obtain
the sections of $L_{\mathbf{m}}$ over $Z_{\mathbf{i}}$. 
We restrict it further down to the dense orbit 
${B \cdot \underline{x}_{\mathbf{i}}}$ of $Z_{\mathbf{i}}$, and then 
by using the orbit map 
\begin{equation*}
B \longrightarrow B \cdot \underline{x} _{\mathbf{i}} \subset Z_{\mathbf{i}}
\end{equation*} 
we pull back the restriction to obtain functions $\xi$
on $B_n =\overline{B}$.

Recall that $x_{i_k}$ in $\underline{x}_{\mathbf{i}}=(x_{i_1},..., x_{i_{\ell}})$ 
is the $|C^{(k)}|$-dimensional subspace of $\mathbb{C}^{n}$ spanned 
by $e_j$ for all $j \in C^{(k)}$. Therefore, the functions $\xi$  
derived from $\delta_j^{(k)}$ are the minors defined on $B_n$ 
with the columns specified by the column set $C^{(k)}$.
This shows that $H^{0}(Z_{\mathbf{i}},L_{\mathbf{m}})$ is spanned 
by tableaux of shape $(\mathbf{i}, \mathbf{m})$ given in \eqref{ttt}.
\end{proof}

Then, we can consider the section ring 
$\mathcal{R}_{\mathbf{i},\mathbf{m}}$ with respect to $L_{\mathbf{m}}$ as the 
$\mathbb{Z}_{\geq 0}$ graded algebra generated by tableaux of shape 
$(\mathbf{i},\mathbf{m})$:%
\begin{equation*}
\mathcal{R}_{\mathbf{i,m}}=\bigoplus_{d\geq 0}\mathsf{M}(\mathbf{i},d\mathbf{m})
\end{equation*}%
where $d\mathbf{m}=(dm_{1},\cdots ,dm_{{\ell}})$. 
We remark that the multiplicative structure of this ring can be 
described by the \textit{straightening laws}, which are in
our case essentially Grosshans-Rota-Stein syzygies given in \cite{DRS76}.
We refer the readers to \cite{Ta01} for more details in this direction.

\medskip


\section{Flat Deformations of the Section Rings}\label{S-flat}


In this section, we describe $\mathbb{C}$-bases of the section spaces and then
prove that the section ring is a flat deformation of a semigroup ring.

\subsection{Straight tableaux}

For a Young diagram $\lambda$, it is well known that semistandard tableaux form
a $\mathbb{C}$-basis of the space spanned by tableaux of shape $\lambda$ (e.g., \cite{DRS76, MS05}).
Now we discuss an analogous result for row convex shape, which is given in \cite{Ta01}
in a general setting of polynomial superalgebras.

\begin{definition}
A row standard \textit{tableau} $\mathsf{t}$ of shape $(\mathbf{i},\mathbf{m})$
given in \eqref{ttt} is called a straight tableau, if it
satisfies the following condition: for two cells $(i,k)$ and $(j,k)$ with 
$i>j$ in the same column, the entry in the upper cell $(i,k)$ may be
strictly larger than the entry in the lower cell $(j,k)$ only if the cell $(i,k-1)$
exists and contains an entry weakly larger than the one in the cell $(j,k)$.
\end{definition}

For example, each of the first three tableaux below can be a part of 
a straight tableau while the last one can not be, because in the last tableau 
$3$ in the second column is less than $4$ in the same column and $1$ left 
to the $4$ is less than $3$:%
\begin{equation*}
\young(\ \ 12,3456,\ \ 5\ ,\ 57\ ),\ \young(\ \ 12,3456,\ \ 3\ ,\ 67\ ),\ %
\young(\ \ 25,3457,\ \ 6\ ,\ 38\ ),\ \young(\ \ 25,1457,\ \ 7\ ,\
38\ )
\end{equation*}

\medskip

A monomial order on the polynomial ring $\mathbb{C}[M_{n}]$ is
called a \textit{diagonal term order} if the leading monomial of a
determinant of any minor defined on $M_{n}$ is equal to the product of the
diagonal elements. For a subring $\mathcal{R}$ of the polynomial ring we let 
$in(\mathcal{R})$ denote the algebra generated by the leading monomials $in(f)$
of all $f\in \mathcal{R}$ with respect to a given monomial order. Note that
the collection of leading monomials forms a semigroup, therefore 
$in(\mathcal{R})$ is a semigroup algebra and $Spec(in(\mathcal{R}))$
is an affine toric variety in the sense of \cite{St95}. Recall that for a subring 
$\mathcal{R}$ of a polynomial ring, a set $\{f_i: i \in I \}$ of 
elements of $\mathcal{R}$ is called a SAGBI basis, if $\{in(f_i): i\in I \}$ 
generates the associated semigroup algebra $in(\mathcal{R})$.

\begin{proposition}\label{STM basis}
Let $D$ be a row-convex shape. 
\begin{enumerate}
\item \cite[Theorem 6.2]{Ta01} Straight tableaux of shape $D$ form 
a $\mathbb{C}$-basis for the space spanned by all the tableaux of shape $D$.
\item \cite[Theorem 7.8]{Ta01} Straight tableaux of shape $D$ form 
a SAGBI basis of the graded algebra $\mathcal{R} \subset \mathbb{C}[M_{n}]$ 
generated by all the tableaux of shape $D$ with respect to 
any diagonal term order.
\end{enumerate}
\end{proposition}

From the fact that the shape $(\mathbf{i},\mathbf{m})$ is row convex, it 
follows from the above Proposition that straight tableaux form a $\mathbb{C}$-basis
of the section ring $\mathcal{R}_{\mathbf{i,m}}$, and that the straight tableaux of shape
$(\mathbf{i},\mathbf{m})$ form a SAGBI basis for
$\mathcal{R}_{\mathbf{i,m}}$. We will study more properties of straight 
tableaux in Section \ref{S-straight}.

\subsection{Flat deformation}

Now we study a flat deformation of the section ring 
$\mathcal{R}_{\mathbf{i}, \mathbf{m}}$. The technique is basically the same as the one for
the Grassmannians and the flag varieties given in, for example, \cite{KM05,St95,MS05}.

\begin{theorem} \label{deformation}
The section ring $\mathcal{R}_{\mathbf{i}, \mathbf{m}}$ of 
the Bott-Samelson variety $Z_{\mathbf{i}}$ can be flatly
deformed into an affine semigroup ring.
\end{theorem}
\begin{proof}
We show that there is a flat $\mathbb{C}[t]$ module 
$\mathcal{R}_{\mathbf{i},\mathbf{m}}^{t}$ whose general fiber is isomorphic to
$\mathcal{R}_{\mathbf{i},\mathbf{m}}$ and special fiber is isomorphic to the
semigroup ring $in(\mathcal{R}_{\mathbf{i},\mathbf{m}})$. 
Lemma \ref{row-convex} shows that any tableau of shape $(\mathbf{i}, \mathbf{m})$ 
with the column sets $\{C_{1}^{(1)},\cdots ,C_{{\ell}}^{({\ell})}\}$ is 
a row-convex tableau. Therefore, we can apply 
Proposition \ref{STM basis} to $\mathcal{R}_{\mathbf{i},\mathbf{m}}$
to conclude that the set of straight tableaux of shape $(\mathbf{i},\mathbf{m})$ 
forms a SAGBI basis for the ring $\mathcal{R}_{\mathbf{i},\mathbf{m}}$
with respect to a diagonal term order. Then, from the existence of a finite SAGBI
basis, by \cite{CHV96}, there exists a $\mathbb{Z}_{\geq 0}$ filtration 
$\{F_{\alpha }\}$ on $\mathcal{R}_{\mathbf{i,m}}$ such that the associated
graded ring of the Rees algebra $\mathcal{R}_{\mathbf{i},\mathbf{m}}^{t}$ with
respect to $\{F_{\alpha }\}$:
\begin{equation*}
\mathcal{R}_{\mathbf{i,m}}^{t}=\bigoplus_{\alpha \geq 0}
F_{\alpha }(\mathcal{R}_{\mathbf{i},\mathbf{m}})t^{\alpha}
\end{equation*}%
is isomorphic to $in(\mathcal{R}_{\mathbf{i,m}})$. Then, by the general
property of the Rees algebra, $\mathcal{R}_{\mathbf{i,m}}^{t}$ is flat over 
$\mathbb{C}[t]$ with general fiber isomorphic to $\mathcal{R}_{\mathbf{i,m}}$ 
and the special fiber isomorphic to the associated graded ring
which is $in(\mathcal{R}_{\mathbf{i,m}})$.
\end{proof}

\medskip


\section{Straight Tableaux and the Space of Sections}\label{S-straight}


In this section, we study more details on the $\mathbb{C}$-basis of the space 
$\mathsf{M}(\mathbf{i}, \mathbf{m})\cong H^{0}(Z_{\mathbf{i}},L_{\mathbf{m}})$ 
given by straight tableaux in Proposition \ref{STM basis}, and then its connection 
to the natural map from the Bott-Samelson variety to the flag variety.

\subsection{Contra-tableaux}

To simplify our notation, we shall keep using the notation 
\begin{equation*}
{\ell}=n(n-1)/2  \text{\ \ and \ \ } p_{j}=j(j-1)/2
\end{equation*} 
for $2\leq j\leq n-1$. Also, fix an arbitrary multiplicity 
$\mathbf{m}=(m_{1},\cdots ,m_{{\ell}})\in \mathbb{Z}_{\geq 0}^{{\ell}}$.

\begin{definition}
A contra-tableau is a filling of a skew Young diagram 
\begin{equation*}
(k,k,\cdots,k)/(\lambda _{1},\lambda _{2},\cdots )
\end{equation*}
with $k\geq \lambda_{1}\geq \lambda _{2}\geq \cdots \geq 0$ such that the entries 
in each column are weakly increasing from top to bottom and the entries in each 
row are strictly increasing from left to right.
\end{definition}

For example, a contra-tableau of shape 
$(4,4,4,4,4)/ (3,3,3,2,1)$ can be encoded in a rectangular array as follows 
\begin{equation*}
\young(\ \ \ 1,\ \ \ 1,\ \ \ 2,\ \ 12,\ 134).
\end{equation*}%
Recall that the usual semistandard tableaux can encode weight basis
elements for irreducible polynomial representations of the general
linear group. Similarly, one can use contra-tableaux to encode
weight vectors of a contragradient representation of an irreducible
polynomial representation of the general linear group. Here, our
goal is to decompose a straight tableau into contra-tableaux.

First, we can decompose the shape $(\mathbf{i}, \mathbf{m})$ into
skew Young diagrams as follows. For $1\leq j\leq n-1$,
let us set $\mathbf{m}(j)=(m_{1}^{\prime },\cdots ,m_{{\ell}}^{\prime })$ where $%
m_{i}^{\prime }=m_{i}$ for $p_{j}<i\leq p_{j+1}$ and $m_{i}^{\prime }=0$
otherwise. Then $\mathbf{m}=\mathbf{m}(1)+\cdots +\mathbf{m}(n-1)$ in $\mathbb{Z}^{\ell}$.

\begin{example}
If $n=4$ and $\mathbf{m}=(1,1,\cdots ,1)\in \mathbb{Z}_{\geq 0}^{6}$, then 
$(\mathbf{i}, \mathbf{m}(1))$, $(\mathbf{i}, \mathbf{m}(2))$, 
$(\mathbf{i}, \mathbf{m}(3))$ respectively correspond to the shapes:%
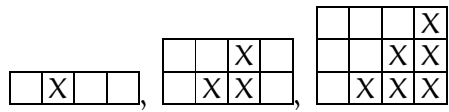
\begin{equation*}
\young(\ X\ \ ),\ \young(\ \ X\ ,\ XX\ ),\ \young(\ \ \ X,\ \ XX,\ XXX)
\end{equation*}%
Note that this is equivalent to the decomposition of the shape 
$(\mathbf{i}, \mathbf{m})$ given in Example \ref{exii} into maximal 
possible Young diagrams.
\end{example}

If $\mathbf{m}=(1,1,\cdots ,1)$, then from the second statement of Lemma \ref%
{row-convex}, the shape $\left(\mathbf{i}, \mathbf{m}(j)\right)$ is a skew
Young diagram $(j+1,j+1,\cdots ,j+1)/(j,j-1,\cdots ,1)$ of length 
$j$. By repeating the $k$-th rows $m_{p_{j}+k}$ times, we have a skew
Young diagram of length $|\mathbf{m}(j)|$. Then from the definition
of straight tableaux, it is straightforward to check that every
straight tableau in a skew diagram is a contra-tableau. See also
\cite[Proposition 4.3]{Ta01}.

\begin{lemma}\label{straight equal contra}
For each $j$, every straight tableau of shape 
$\left( \mathbf{i}, \mathbf{m}(j)\right)$ is a contra-tableau.
\end{lemma}

Note that this lemma shows that
the basis of the space $\mathsf{M}(\mathbf{i}, \mathbf{m}(j)) \cong 
H^{0}(Z_{\mathbf{i}},L_{\mathbf{m}(j)})$ is simply given by contra-tableaux 
tableaux, and then as a consequence we can obtain a description of elements 
in the section ring $\mathcal{R}_{\mathbf{i,m}}$ as products of
contra-tableaux. That is, we have a natural projection
\begin{equation}\label{projection}
\mathsf{M}(\mathbf{i}, \mathbf{m}(1))\otimes \cdots \otimes 
\mathsf{M}(\mathbf{i}, \mathbf{m}(n-1)) \rightarrow 
\mathsf{M}(\mathbf{i}, \mathbf{m})  
\end{equation}%
sending $\mathsf{t}_{1}\otimes \cdots \otimes \mathsf{t}_{n-1}$ 
to the product $\mathsf{t}_{1}\cdot ...\cdot \mathsf{t}_{n-1}\in 
\mathsf{M}(\mathbf{i}, \mathbf{m})$ where $\mathsf{t}_{j}$ is 
a contra-tableau in $\mathsf{M}(\mathbf{i}, \mathbf{m}(j))$ for each $j$.

For example, if $n=4$ and $\mathbf{m}=(1,2,1,1,1,3)$, then the
product map $\mathsf{t}_{1}\otimes \mathsf{t}_{2} \otimes
\mathsf{t}_{3} \rightarrow \mathsf{t}$ gives
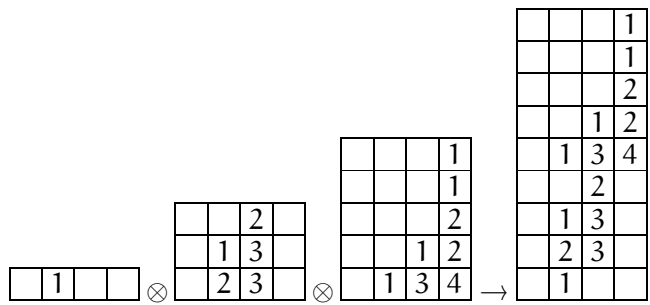
\begin{equation*}
\young(\ 1\ \ )\otimes \young(\ \ 2\ ,\ 13\ ,\ 23\ )\otimes \young(\
\ \ 1,\ \ \ 1,\ \ \ 2,\ \ 12,\ 134)\rightarrow \young(\ \ \ 1,\ \ \
1,\ \ \ 2,\ \ 12,\ 134,\ \ 2\ ,\ 13\ ,\ 23\ ,\ 1\ \ )
\end{equation*}%
Note that the product is not a straight tableau, but it can be expressed by
a linear combination of straight tableaux in 
$\mathsf{M}(\mathbf{i}, \mathbf{m})\subset \mathcal{R}_{\mathbf{i,m}}$ 
by successive application of straightening laws
mentioned after Proposition \ref{coordinate ring}.

\medskip

\subsection{Projection to $G/B$}

Now we discuss the natural map from the Bott-Samelson variety 
$Z_{\mathbf{i}}$ to the flag variety $G/B$ in terms of our basis description.

The projection map \eqref{projection} is compatible with the decomposition
of a straight tableau of shape $(\mathbf{i}, \mathbf{m})$ into
contra-tableaux. More precisely, a straight tableau $\mathsf{t}$ of shape 
$\left( \mathbf{i}, \mathbf{m}\right)$ can be factored into a product 
$\mathsf{t}_{1}\cdot ...\cdot \mathsf{t}_{n-1}$ of straight tableaux 
$\mathsf{t}_{j}$\ of shape $(\mathbf{i}, \mathbf{m}(j))$ for $1\leq j\leq n-1$. 
Then, by Lemma \ref{straight equal contra}, $\mathsf{t}_{j}$ are 
contra-tableaux for all $j$.

In particular, for each $1\leq j\leq n-2$, let us consider a straight
tableau $\mathsf{t}_{j}^{0}$ of shape $\left( \mathbf{i}, \mathbf{m}(j)\right)$ 
such that for each $a$ and $b$ such that $1\leq b\leq m_{j}$ and 
$p_{j}+1\leq a\leq p_{j+1}$, the row indices and the column indices are
equal: $R_{b}^{(a)}=C^{(a)}$, i.e.,%
\begin{eqnarray*}
\mathsf{t}_{1}^{0}&=&[C^{(2)}:C^{(2)}]^{m_{1}}; \\
\mathsf{t}_{j}^{0}&=&[C^{(p_{j}+1)}:C^{(p_{j}+1)}]^{m_{p_{j}+1}}\cdot
[C^{(p_{j}+2)}:C^{(p_{j}+2)}]^{m_{p_{j}+2}}\cdot ...\cdot
[C^{(p_{j+1})}:C^{(p_{j+1})}]^{m_{p_{j+1}}}
\end{eqnarray*}%
for $2 \leq j \leq n-2$. This is equivalent to say that
$\mathsf{t}_{j}^{0}$ is obtained by filling
in all the cells corresponding to the subshapes 
$\left( \mathbf{i}, \mathbf{m}(j) \right)$ of the shape 
$\left( \mathbf{i}, \mathbf{m}\right)$ with maximum possible numbers.

Then, for any contra-tableau $\mathsf{t}$ of shape $(\mathbf{i}, \mathbf{m}(n-1))$, 
we can find a straight tableau $\widehat{\mathsf{t}}$ of shape 
$(\mathbf{i}, \mathbf{m})$ such that%
\begin{equation*}
\widehat{\mathsf{t}}=\left( \mathsf{t}_{1}^{0}\cdot ...\cdot \mathsf{t}%
_{n-2}^{0}\right) \cdot \mathsf{t}
\end{equation*}%
and this provides the following injection:%
\begin{eqnarray}\label{injection}
H^{0}(G/B,L_{\mathbf{\lambda }}) &\rightarrow &\mathsf{M}(\mathbf{i}, \mathbf{m}) \\
\mathsf{t} &\mapsto &\left( \mathsf{t}_{1}^{0}\cdot ...\cdot \mathsf{t}%
_{n-2}^{0}\right) \cdot \mathsf{t} \notag
\end{eqnarray}%
where $H^{0}(G/B,L_{\mathbf{\lambda }})$ is the space of section of the line
bundle $L_{\mathbf{\lambda}}$ on $G/B$ and $\lambda $ is the dominant
weight determined by the shape $\mathbf{m}(n-1)$ as a Young diagram.

For example,
\begin{equation*}
\young(\ \ \ 1,\ \ \ 1,\ \ \ 2,\ \ 12,\ 134)\rightarrow 
\young(\ \ \ 1,\ \ \ 1,\ \ \ 2,\ \ 12,\ 134,\ \ 3\ ,\ 23\ ,\ 23\ ,\ 2\ \ )
\end{equation*}

Finally, by extending the map \eqref{injection}, we have

\begin{proposition} There is a natural map from the section ring of the 
flag variety to the section ring of $Z_{\mathbf{i}}$
\begin{equation*}
\bigoplus_{d \geq 0} H^{0}(G/B,L^d_{\mathbf{\lambda }}) \longrightarrow 
\mathcal{R}_{\mathbf{i,m}}=\bigoplus_{d\geq 0}H^{0}
(Z_{\mathbf{i}},L_{\mathbf{m}}^{d}).
\end{equation*}

\end{proposition}


\section{Three-dimensional Example: Toric Degenerations}\label{S-examp}


In this section we will consider explicit examples of toric
degenerations of a three-dimensional Bott-Samelson variety, and compute
the corresponding Hilbert polynomials. 

\subsection{}

Let $P_1$ and $P_2$ be the following parabolic subgroups 
of $\rm{GL}_{3}(\mathbb{C})$:
\begin{equation*}
P_1=\left(
\begin{array}{ccc} 
* & * & * \\ * & * & * \\ 0 & 0 & *
\end{array} \right), 
\ \ \ 
P_2=\left(\begin{array}{ccc}
 * & * & * \\ 0 & * & * \\ 0 & * & * 
\end{array} \right)
\end{equation*}
and denote by ${\bar P}_1$ and ${\bar P}_2$ their closures in
the space $M_3$ of $3\times 3$ matrices. 

Let $Z$ be the Bott-Samelson variety defined as in \S \ref{Bquot} with $n=3$ 
and ${\mathbf{i}}=(121)$. That is, 
\begin{equation*}
Z = P_1\times P_2 \times P_1 / B^3
\end{equation*}
with the action of  $B^3$:
\begin{equation*}
(p_1, p_2, p_3) \cdot (b_1, b_2, b_3) = (p_1 b_1, b_1^{-1}p_2 b_2, b_2^{-1}p_3 b_3).
\end{equation*}
It can also be viewed as an invariant theory quotient
of the product of the closures ${\bar P}_1\times {\bar P}_2 \times
{\bar P}_1$ by the action of $B^3$ in the obvious way.

We will denote the elements of the first copy of $P_1$ by
\begin{equation*}
p_1=\left(\begin{array}{ccc} a_{11} & a_{12} & a_{13} \\
a_{21} & a_{22} & a_{23} \\ 0 & 0 & a_{33} \end{array} \right),
\end{equation*}
the elements of $P_2$ by
\begin{equation*}
p_2=\left(\begin{array}{ccc} b_{11} & b_{12} & b_{13} \\
0 & b_{22} & b_{23} \\ 0 & b_{32} & b_{33} \end{array} \right),
\end{equation*}
and the elements of the second copy of $P_1$ by
\begin{equation*}
p_3=\left(\begin{array}{ccc} c_{11} & c_{12} & c_{13} \\
c_{21} & c_{22} & c_{23} \\ 0 & 0 & c_{33} \end{array} \right).
\end{equation*}
The same notation will be used for the elements of their closures 
in $M_3$.

\subsection{}

Next, we will describe a Pl\"ucker-type embedding of $Z$ into the
product of three projective spaces:
\begin{equation*}
{\mathcal H}:={\rm Proj}(s_1, s_2)\times {\rm Proj}(r_{23}, r_{13},
r_{12})\times {\rm Proj}(q_1, q_2, q_3)\simeq {\mathbb C}{\mathbb P}^1\times
 {\mathbb C}{\mathbb P}^2\times {\mathbb C}{\mathbb P}^2.
\end{equation*}

\smallskip

Let a point in $Z$ be represented by three matrices $(p_1, p_2,
p_3)$ in the above form, then we denote by 
$s_i$ the $1\times 1$ minor of the matrix $p_1$ with column 1 and row $i$. 
Therefore, we have 
\begin{equation*}
s_1=a_{11}  \text{\ \ and \ } s_2=a_{21}. 
\end{equation*}
($s_3$ would be identically equal to zero, so we do not use it.)
Next, we denote by $r_{ij}$ the $2\times 2$ minor of the matrix $p_1p_2$ with
columns $1,2$ and rows $i,j$. Explicitly,
\begin{align*}
&r_{12} =a_{11}b_{11}(a_{22}b_{22} + a_{23}b_{32}) - a_{21}b_{11}(a_{12}b_{22}+a_{13}b_{32}),\\
&r_{13} =a_{11}a_{33}b_{11}b_{32}, \text{\ \ and\ \ } r_{23} =a_{21}a_{33}b_{11}b_{32}.
\end{align*}
Finally, we denote by $q_{i}$ the $1\times 1$ minor of the matrix 
$p_1p_2p_3$ with column $1$ and row $i$:
\begin{align*}
&q_1 =
a_{11}b_{11}c_{11}+(a_{11}b_{12}+a_{12}b_{22}+a_{13}b_{32})c_{21},\\
&q_2 =
a_{21}b_{11}c_{11}+(a_{21}b_{12}+a_{22}b_{22}+a_{23}b_{32})c_{21}, \\
&\text{and\ \ } q_3 = a_{33}b_{32}c_{21}.
\end{align*}

Then, $Z$ can be viewed as a subvariety of ${\mathcal H}$, the
product of three projective spaces, defined by the following two
homogeneous equations (or Pl\"ucker relations):
\begin{equation}\label{2Plucker}
s_1r_{23} - s_2r_{13}=0 \text{\ \ and\ \ } q_1r_{23} - q_2r_{13}+q_3r_{12}=0.
\end{equation}

\begin{proposition}\label{propHP}
The Hilbert polynomial of $Z$ is given by
\begin{equation*}
{\rm HP}_Z(s) = \frac{5s^3+11s^2+8s+2}{2}.
\end{equation*}
\end{proposition}
\begin{proof}
Let, as before,
${\mathcal H}={\mathbb C}{\mathbb P}^1\times {\mathbb C}{\mathbb P}^2
\times {\mathbb C}{\mathbb P}^2$ and let $\pi_1$, $\pi_2$ and $\pi_3$ stand
for the projections onto the corresponding factors. Write
\begin{equation*}
L=\pi_1^*({\mathcal O}(1)), \  M_1=\pi_2^*({\mathcal O}(1)), 
\text{\ \ and\ \ } M_2=\pi_3^*({\mathcal O}(1)). 
\end{equation*}
We will also denote by the same letters $L$, $M_1$, and $M_2$ the corresponding 
classes of divisors in the Chow ring of ${\mathcal H}$. Let $X$ be the element 
of the Chow ring of ${\mathcal H}$ corresponding to $Z$, and let 
\begin{equation*}
D = n(L+M_1+M_2). 
\end{equation*}
For large enough integral values of $s$, the Hilbert polynomial ${\rm HP}_Z(s)$ 
coincides with ${\rm dim}(H^0(sD_{\vert Z}))$, which, due to vanishing,
is the same as the Euler characteristic of $sD_{\vert Z}$.

The Riemann-Roch theorem for smooth Fano threefolds (e.g., \cite{IP99}) 
asserts that 
\begin{equation*}
\chi(nD_{\vert Z})=\frac{D^3_{\vert Z}}{6}n^3-\frac{D^2_{\vert Z}K_Z}{4}n^2
+ \frac{D_{\vert Z}(K_Z^2+c_2(Z))}{12}n+1.
\end{equation*}
Now, $X=(L+M_1)(M_1+M_2)$, therefore, by the adjunction formula we get
$-K_Z=(L+M_1+2M_2)_{\vert {Z}}$ and hence $(L+M_1+2M_2)_{\vert Z}c_2(Z)=24$.
To find $(M_2)_{\vert Z}c_2(Z)$, we will use the same Riemann-Roch formula, 
but for $M_2$, note that ${\rm dim}(H^0(M_2))=3$.
Finally, the intersection products satisfy 
\begin{equation*}
L^2=0, \ \ M_1^3 = M_2^3=0, \text{\ \ and \ \ } LM_1^2M_2^2=1, 
\end{equation*}
which leads to a straightforward computation of the required polynomial.
\end{proof}

\subsection{}

Forgetting the first component, one can consider the projection:
\begin{equation*}
{\mathcal H}\to {\mathbb C}{\mathbb P}^2\times {\mathbb C}
{\mathbb P}^2.
\end{equation*}
The image of $Z$ under this projection is naturally 
the $3$-dimensional flag variety ${\rm Fl}_3$, 
sitting inside ${\mathbb C}{\mathbb P}^1\times {\mathbb C}{\mathbb P}^2$ 
as the zero set of the second Pl\"ucker relation in \eqref{2Plucker}.

There are two naive ways to construct toric degenerations of $Z$,
the first is to consider the family of varieties, parameterized by
$\tau\in{\mathbb C}$, where the second equation is modified to
\begin{equation*}
q_1r_{23} - q_2r_{13} + \tau q_3r_{12}=0.
\end{equation*}
One can easily notice that the special toric fiber of this family,
corresponding to $\tau=0$ is a reducible variety and has two
irreducible components: one, denoted by ${\mathcal G}$, is
isomorphic to ${\mathbb C}{\mathbb P}^1\times {\mathbb C}{\mathbb P}^2$, 
and corresponds to $r_{23}=r_{13}=0$, and the second component, denoted by $D_3$, 
a three-dimensional toric variety, which is actually non-singular. 
Combinatorially, the moment polytope for $D_3$ is a cube, and is drawn 
schematically on the figure below.
(To simplify computations, we assumed that the members of the family
are polarized by the invertible sheaf induced from 
${\mathcal O}(1)\times {\mathcal O}(1)\times {\mathcal O}(1)$ 
on ${\mathcal H}$.)

\begin{center}
\begin{picture}(140,43)(0,0)
  \put(50,3){
    \unitlength=7mm
    \drawpolygon(0,0)(0,3)(2,3)(2,1)
    \drawline[AHnb=0](2,1)(5,4)(5,5)(3,5)(0,3)
    \drawline[AHnb=0](2,3)(5,5)
    \drawline[AHnb=0,dash={0.2}0](0,0)(3,3)(3,5)
    \drawline[AHnb=0,dash={0.2}0](3,3)(5,4)
  }
\end{picture}
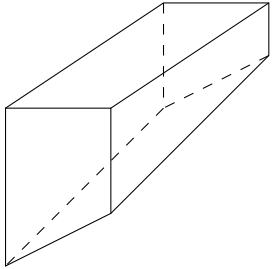
FIGURE 1.
\end{center}
\medskip

The fact that the special fiber of this flat family of varieties
over ${\mathbb C}$ is reducible and is given by the union of two
non-singular components is quite amusing. The intersection of these
two components is a smooth two-dimensional toric variety, denoted by
$K_2$, known as the Hirzebruch surface of degree one.

One can compute the Hilbert polynomials for the chosen polarization,
denoted by ${\rm HP}$ of the irreducible components, which are known
(e.g., \cite{MS05}) to be the same as the Ehrhart polynomials, denoted by
${\rm EP}$ of their moment polytopes,
as well as their Ehrhart series, ${\rm ES}$. 
Using a computer program \cite{Latte}, we have obtained:
\begin{align*}
&{\rm ES}(D_3) = \frac{3t^2+8t+1}{(1-t)^4}, \ \ 
{\rm EP}(D_3)  = {\rm HP}(D_3) = 2s^3+5s^2+4s+1 = (s+1)^2(2s+1),\\
&{\rm ES}({\mathcal G}) = \frac{1+2t}{(1-t)^4}, \ \ \ 
{\rm EP}({\mathcal G}) = {\rm HP}({\mathcal G}) = \frac{s^3+4s^2+5s+2}{2}
= \frac{(s+1)^2(s+2)}{2}, \\
&{\rm ES}(K_2) = \frac{1+2t}{(1-t)^3}, \ \ \ 
{\rm EP}(K_2) = {\rm HP}(K_2) = \frac{3s^2+5s+2}{2} = \frac{(s+1)(3s+2)}{2},
\end{align*}
and this allows us to check that
\begin{align*}
{\rm HP}(Z) &= {\rm HP}(D_3) + {\rm HP}({\mathcal G}) - {\rm HP}(K_2) \\
&= \frac{5s^3+11s^2+8s+2}{2} = \frac{(s+1)(5s^2+6s+2)}{2}.
\end{align*}

This fact was also verified, independently, using a software package
\cite{Sing}, by representing $Z$ as a subvariety in 
${\mathbb C}{\mathbb P}^{17}$ via Segre embedding, defined by the following 95
equations, where $[a_1:\cdots :a_9:b_1:\cdots :b_9]$ are the homogeneous
coordinates on ${\mathbb C}{\mathbb P}^{17}$:
\begin{equation*}
a_ib_j=a_jb_i\ \  {\rm for} \ \ 1\le i < j \le 9,
\end{equation*}
similarly,
\begin{align*}
& a_ka_l=a_ma_n, \ a_kb_l=a_mb_n, \ a_kb_l=b_ma_n, \\ 
& b_ka_l=a_mb_n, \ b_ka_l=b_ma_n, \ b_kb_l=b_mb_n,
\end{align*}
for the following nine choices of quadruples of indexes $(k,l,m,n)$:
\begin{align*}
& (1,5,2,4), \ (1,6,3,4), \ (2,6,3,5,), \ (1,8,2,7), \ (1,9,3,7),\\
& (2,9,3,8), \ (4,8,5,7), \ (4,9,6,7), \ {\rm and} \ (5,9,6,8),
\end{align*}
and the last five:
\begin{align*}
& a_1+b_4=0, \ \ a_2 + b_5 = 0, \ \ a_3 + b_6 = 0, \\
& a_1+a_5+a_9 = 0,  \ {\rm and} \  b_1+b_5+b_9 =0.
\end{align*}

\subsection{}

The second way to obtain a flat toric degeneration of $Z$ is to
consider a different family of varieties inside ${\mathcal H}$, also
parameterized by $\tau\in{\mathbb C}$ and given by the following two
equations:
\begin{equation*}
s_1r_{23}-s_2r_{13}=0 \ \ \ \text{\ and \ } \ \ \ 
q_1r_{23}-\tau q_2r_{13}+q_3r_{12}=0.
\end{equation*}
One can see that the special fiber of this flat family,
corresponding to $\tau=0$, is a singular toric variety, denoted by $Y_3$, 
whose moment polytope combinatorially is represented by the picture drawn below:
\begin{center}
\begin{picture}(140,43)(0,0)
  \put(50,3){
    \unitlength=10mm
    \drawpolygon(0,0)(0,2)(2,2)(2,1)(1,0)
    \drawline[AHnb=0](1,0)(4,2)(5,3)(2,1)
    \drawline[AHnb=0](0,2)(3,3)(5,3)(2,2)
    \drawline[AHnb=0,dash={0.2}0](0,0)(3,2)(3,3)
    \drawline[AHnb=0,dash={0.2}0](3,2)(4,2)
  }
\end{picture}
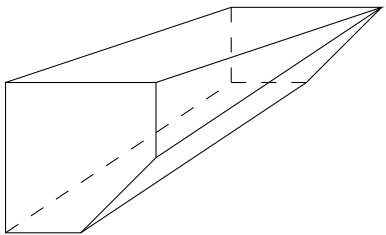
FIGURE 2.
\end{center}

\medskip

The Ehrhart series and the Ehrhart polynomial of the moment polytope
of special fiber corresponding to the same, previously chosen,
polarization, is given by
\begin{equation*}
{\rm ES}(Y_3) = \frac{5t^2+9t+1}{(1-t)^4} \text{\ \ and \ } 
{\rm EP}(Y_3) = \frac{5s^3+11s^2+8s+2}{2}.
\end{equation*}
Not surprisingly, again we see that ${\rm EP}(Y_3) = {\rm HP}(Z)$.

\medskip

\subsection*{Acknowledgment} We thank Mikhail Kogan for his contribution on an early
stage of the project. We also thank Ivan Cheltsov for help with 
Proposition \ref{propHP}.

\end{document}